\RequirePackage{amsmath}
\documentclass[runningheads]{llncs}
\usepackage[T1]{fontenc}
\usepackage{graphicx}
\usepackage{xcolor}

\usepackage[ruled,vlined,linesnumbered]{algorithm2e}
\usepackage{subcaption}
\usepackage{enumitem}
\usepackage{amssymb}
\usepackage[]{hyperref}
\hypersetup{
    pdftitle={Computational Tradeoffs of Optimization-Based Bound Tightening in ReLU Networks},
    pdfauthor={Badilla et al},
    pdfkeywords={Activation Bounds, Neural Networks, Neural Surrogate Models, Mixed-Integer Linear Programming, Rectified Linear Units},
    bookmarksnumbered=true,
    bookmarksopen=true,
    bookmarksopenlevel=2,
    pdfstartview=Fit,
    pdfpagemode=UseOutlines
}

\newcommand{\linf}[1]{
    ||#1||_\infty
}

\begin{document}
\title{Computational Tradeoffs of Optimization-Based Bound Tightening in ReLU Networks}
\titlerunning{Computational Trade-offs of OBBT in ReLU Networks}
\author{Fabián Badilla\inst{1,2} \and 
Marcos Goycoolea\inst{3} \and 
Gonzalo Muñoz\inst{4,5} \and 
Thiago Serra\inst{2}
}
\authorrunning{F. Badilla et al.}
\institute{Universidad de Chile, Santiago, Chile \and
Bucknell University, Lewisburg PA, United States \and
Pontificia Universidad Católica de Chile, Santiago, Chile\and
Universidad de O'Higgins, Rancagua, Chile \and 
Instituto Sistemas Complejos de Ingenier\'ia (ISCI), Chile
}%
\maketitle
\begin{abstract}
The use of Mixed-Integer Linear Programming (MILP) models to represent neural networks with Rectified Linear Unit (ReLU) activations has become increasingly widespread in the last decade.
This has enabled the use of MILP technology to test---or stress---their behavior, to adversarially improve their training, and to embed them in optimization models leveraging their predictive power.
Many of these MILP models rely on \emph{activation bounds}. That is, bounds on the input values of each neuron.
In this work, we explore the tradeoff between the tightness of these bounds and the computational effort of solving the resulting MILP models. 
We provide guidelines for implementing these models based on the impact of network structure, regularization, and rounding.

\end{abstract}

\section{Introduction}

One of the key sources of the expressiveness of neural networks is their nonlinearity, which comes from the activation functions. The most commonly used activation function in the last decade is the Rectified Linear unit (ReLU) activation function~\cite{hahnloser2000origin,nair2010rectified,glorot2011rectifier,lecun2015nature,ramachandran2018pop}, defined as $\sigma(x) = \max\{0,x\}$.
All ReLU neural networks represent piecewise linear functions~\cite{pascanu2013on}. Conversely, all piecewise linear functions can be represented as ReLU networks~\cite{arora2018understanding}. 
Moreover, ReLU networks are universal function approximators under a variety of architectures~\cite{yarotsky2017relu,lu2017expressive,hanin2017approximating,park2021width}, meaning that they can represent almost any function to arbitrary precision in the same way as with other activation functions~\cite{cybenko1989approximation,funahashi1989approximate,hornik1989approximator}. 

This piecewise linear structure of the ReLU activation function facilitates embedding these networks in mixed-integer linear models. For an input $x\in [\ell, u]$, we can represent such $x$--$y$ relation $y=\sigma(x)$ as the following linear system.
\begin{equation}
    \label{single-bigM}
    x = y - h,\, y\leq u z,\, h \leq -\ell (1-z),\, y \geq 0, h \geq 0,\, z\in \{0,1\}.
\end{equation}
Here, $z$ indicates if $x\geq 0$ or not.
This system can model every neuron using ReLU, thus effectively representing the neural network evaluation as a system of linear inequalities over mixed-integer variables. 

These models have been used in several applications involving trained neural networks in one way or another. 
The applications focused on the neural networks themselves include the generation of adversarial inputs~\cite{cheng2017maximum,fischetti2018deep,khalil2018combinatorial}, 
verification and robustness certification~\cite{Tjeng2017Nov,li2022sok,liu2021algorithms}, 
feature visualization~\cite{fischetti2018deep}, 
counting linear regions~\cite{serra2017bounding,serra2020empirical,cai2023pruning}, 
compression~\cite{serra2020lossless,serra2021scaling,elaraby2023oamip}, 
providing strong convex relaxations~\cite{anderson2018strong,anderson2020strong,tsay2021partition}, 
generating heuristic solutions for neural surrogate models~\cite{perakis2022optimizing,tong2023walk}, 
producing counterfactual explanations~\cite{kanamori2021counterfactual}, 
and constraining reinforcement learning policies~\cite{delarue2020rlvrp,burtea2023safe}. 
We refer the reader to the recent survey \cite{huchette2023deep} for a review of further uses of polyhedral theory and MILP technology in the specific context of deep learning. 
Moreover, modeling frameworks~\cite{bergman2022janos,ceccon2022omlt,fajemisin2023ocl,maragno2023mixed,gurobi2023ml} have been inspired by the wider variety of problems beyond deep learning, e.g., in power grids~\cite{chen2020voltage,murzakhanov2022powerflow}, 
automated control~\cite{say2017planning,wu2020scalable,yang2021control}, 
vehicle routing~\cite{delarue2020rlvrp},
among others.

An important aspect of these models is that they require valid bounds on the input of each ReLU function, as in \eqref{single-bigM}. 
While the tightness of these bounds does not affect the correctness of the formulation, 
exceedingly large bounds may negatively affect the computational performance when solving the resulting MILP models. This is because the effectiveness of MILP technologies largely depends on the strength of LP relaxation bounds. 

In this work, we study the tradeoffs between the time spent computing tight activation bounds and the time gained from using such bounds in mixed-integer linear models involving trained neural networks.
First, we analyze how much more expensive it is to compute tight bounds using an Optimization-Based Bound Tightening (OBBT) approach compared to a simpler interval-based propagation approach. 
We evaluate runtimes according to different characteristics of the neural networks, such as architectural structure, regularization during training, and network pruning level. 
Second, we test these different bounds when applied to network verification, 
hence characterizing the aforementioned tradeoff.
Our goal is to provide computational guidelines for researchers and practitioners when using MILP models involving trained neural networks.

Though many papers tend to use MILP optimization for bound tightening ~\cite{fischetti2018deep,serra2020empirical,tsay2021partition,Cacciola2023Jul}, to our knowledge there is very little work analyzing this computational tradeoff. An exception is the work of Liu et al.~\cite{liu2021algorithms}, which compares various algorithms for computing activation bounds, some of which we also consider here. 
In this paper, we complement the work of \cite{liu2021algorithms} in two important ways: (i) we include a sensitivity analysis of different bound computation approaches with respect to multiple network characteristics; and (ii) we analyze the effect of the bounds quality in other optimization problems.

\section{Background}

\subsection{Neural Networks}
A neural network is a function \(f\) defined over a
directed graph that maps inputs \(x \in \mathbb{R}^{n_0}\) to \(f(x) \in
\mathbb{R}^{n_L}\). The directed graph \(G = (V,E)\) representing the network decomposes into layers $V_i$ with $i\in [L]_0:=\{0,...,L\}$ such that \(V = \bigcup\limits_{i \in [L]_0} V_i\); \(V_0\) is the input layer and \(V_L\) as the output layer.
Layers $1$ through $L-1$ are the \emph{hidden layers}. 
In this work, we consider fully connected feed-forward networks, meaning that all nodes in layer $V_j$ are connected to all nodes in layer $V_{j+1}$ and that arcs only connect nodes in adjacent layers. This is mainly to simplify the discussion, as this framework can also be adapted to other feed-forward architectures.

The neural network computes an output in the following way. Each layer $l$ outputs a vector $h^{(l)}$ to the following layer, starting with $h^{(0)} = x$. Then we apply an affine transformation ${W^{(l)}}^\top h^{(l)} + b^{(l)} =: a^{(l+1)}$. Finally, the activation function is applied (component-wise) to $a^{(l+1)}$ to produce the output of layer $l+1$, i.e., $h^{(l+1)} = \sigma(a^{(l+1)})$.
The last layer often applies an extra transformation, such as $\mbox{softmax}(z) = e^{z_i}/\sum_{j=1}^w e^{z_j}$, to model a probability distribution. The softmax layer cannot be encoded using MILP, but we can still determine the class by comparing the inputs of softmax since it is a monotonic function.

\subsection{Computing Activation Bounds}
If the input $x\in \mathbb X$ is bounded, it is easy to see that each value of $a^{(l)}$ is bounded as well; these \emph{activation bounds} are needed for a formulation using \eqref{single-bigM}.
Fix a layer $m$ and suppose we have computed values $\mbox{LB}_{j}^{(l)}, \mbox{UB}_{j}^{(l)}$ such that $a_{j}^{(l)}\in [\mbox{LB}_{j}^{(l)}, \mbox{UB}_{j}^{(l)}]$ for $l\leq m-1$.
Then we can compute a valid bound for neuron $n$ as follows:
\begin{subequations}\label{eq:OBBT}
\begin{align}
  \mbox{UB}_{n}^{(m)} =  \max_{x\in\mathbb{X}} \quad & a^{(m)}_n\\
    \text{s.t.} \quad & {W^{(l)}}^\top h^{(l)} + b^{(l)} = a^{(l+1)} && \forall l\in[m-1]_0  \label{cstr:pv_fc}\\
    &  a^{(l)} = h^{(l)}-\overline{h}^{(l)}  && \forall l\in[m-1] \label{cstr:pvnegpos}\\
    &  h_j^{(l)} \leq \mbox{UB}_j^{(l)} z_j^{(l)}  && \forall l\in[m-1], \forall j\in [n_l] \label{cstr:pvub}\\
    & \overline{h}_j^{(l)} \leq -\mbox{LB}_j^{(l)} (1-z_j^{(l)})  && \forall l\in[m-1], \forall j\in [n_l] \label{cstr:pvlb}\\
    & h^{(l)}, \overline{h}^{(l)} \geq 0,\, z^{(l)} \in \{0,1\}^{n_l} && \forall l\in [m-1] \label{cstr:pvrelu}\\
    & h_j^{(0)}\in [\mbox{LB}_j^{(0)}, \mbox{UB}_j^{(0)}] && \forall j \in [n_0] &\label{cstr:pv_initial}
\end{align}
\end{subequations}
Here $n_l := |V_l|$ is the number of neurons in layer $l$ (the width). 
Note also that $\mbox{LB}_{n}^{(m)}$ can be computed by minimizing instead of maximizing. 
Solving \eqref{eq:OBBT} to determine $\mbox{LB}_{j}^{(l)}, \mbox{UB}_{j}^{(l)}$ is NP-hard \cite{katz2017reluplex}, and so it is unlikely that an efficient algorithm exists for solving \eqref{eq:OBBT}. However, solving the linear relaxation of \eqref{eq:OBBT} also provides valid, quickly computable, but weaker bounds. Henceforth, we refer to these weaker bounds as $\widehat{\mbox{LB}}_{j}^{(l)}$ and $\widehat{\mbox{UB}}_{j}^{(l)}$.

Note that even weaker bounds can be computed, without the need to solve any optimization problem. This follows from Proposition \ref{prop:naive}, where upper activation bounds for a node can be computed from those in the previous layer. Also, by reversing this inequality we have lower activation bounds. If we do this for all the layers iteratively, we will get valid bounds for the entire network.

\begin{proposition}[Naive Bound]\label{prop:naive}
    Given $\mbox{UB}^{(l)}$ valid bounds for layer $l\geq 0$, then the following is a valid activation bound for every neuron on layer $l+1$:
$$\linf{a^{(l+1)}} \leq \linf{W^{(l)}}\max_{j\in [n_l]} \left(\mbox{UB}_j^{(l)} \right) + \linf{b^{(l)}}$$
\end{proposition}

To find all activation bounds using formulation \eqref{eq:OBBT}, we proceed forward from the input layer to compute upper and lower activation bounds one layer at a time. This procedure is described in Algorithm \ref{alg:bounder}. We say that the bounds $\mbox{LB}_{n}^{(m)}$ and $\mbox{UB}_{n}^{(m)}$ resulting from this algorithm are the \emph{strong} node activation bounds. Note that in step \ref{step:solve} we can solve the relaxed version of \eqref{eq:OBBT} if we want to compute $\widehat{\mbox{LB}}_{j}^{(l)}$ and $\widehat{\mbox{UB}}_{j}^{(l)}$. We call these the \emph{weak} node activation bounds. If, instead, we ignore step \ref{step:solve} altogether, we obtain what we call the \emph{naive} node activation bounds.

We remark that in our implementation of Algorithm~\ref{alg:bounder}, we parallelized the for loop in line 5; this is possible because problem \eqref{eq:OBBT} is independent $\forall n\in [n_m]$ for a fixed layer $m\in [L]$. Note also that if the MILP solver fails to produce an optimal integer solution, we ask the solver to provide us with the best dual bound available, if there is one, before using the bound from Proposition \ref{prop:naive}. This could happen in rare cases, for example, if the time limit is reached without a dual bound being reported or if there are numerical issues, among others.

This algorithm could be extended to recurrent neural networks with a finite number of steps, but for the purposes of this paper, we will restrict the scope to feedforward neural networks.

\SetKwInput{KwInput}{Input}                % Set the Input
\SetKwInput{KwOutput}{Output}              % set the Output
\begin{algorithm}[t]
\DontPrintSemicolon
\caption{Bounder}\label{alg:bounder}
\KwInput{Trained neural network, input layer bounds $\mbox{LB}^{(0)}, \mbox{UB}^{(0)}$}
Get admissible bounds $\mbox{LB}^{(l)}, \mbox{UB}^{(l)}$ $\forall l\in [L]$ using Proposition \ref{prop:naive}\\
Initialize the variables of \eqref{eq:OBBT} and add \eqref{cstr:pv_initial}\\
\For{$1 \leq m \leq L$}{
    Add constraints $\eqref{cstr:pv_fc}-\eqref{cstr:pvlb}$ for the layer $m$\\
    \For{$1 \leq n \leq n_m$\label{bounder:nodesfor}}{
        Solve \eqref{eq:OBBT} and its minimization versions for the current node $(n,m)$\\ \label{step:solve}
        \If{Valid bounds for $a_n^{(m)}$ are found}{
            Set $\mbox{UB}_n^{(m)}$ and $\mbox{LB}_n^{(m)}$ with the corresponding values\\
        }
        \Else{
            Propagate worst bound of layer $m-1$ with Proposition \ref{prop:naive}.
        }
    }
}
\textbf{Return} $\left\{\mbox{LB}^{(l)}, \mbox{UB}^{(l)}\right\}_{l=0}^{L}$
\end{algorithm}

\subsection{Verification Problems}\label{sec:verification}

To analyze and compare the impact of the different ways of computing the activation bounds, we will use them in a verification problem formulation.
In this problem, we start with reference input $x_0$ and formulate an optimization problem that finds a small perturbation that is classified differently to $x_0$.

For the sake of brevity, we will not display the full model, but it suffices to say that we simply modify the formulation \eqref{eq:OBBT} to have an objective function that maximizes the difference in the classification of $x_0$ and the perturbed input. 

%==============================================================================
%==============================================================================

\section{Experimental Setup}

Our computational experiment consisted of three steps. First, we trained 96 neural networks using the MNIST dataset \cite{mnist} as follows:
\begin{enumerate}
    \item \emph{Architecture:} We used four different neural network architectures, varying in the number of hidden layers \(N_l\) and nodes per layer \(N_n\), represented as \(N_l \times N_n\); each layer being fully-connected and the last being a softmax layer.
    \item \emph{Regularization and Training Duration:} The networks were trained for 50 epochs with L1 and L2 regularization, employing a penalty parameter \(\lambda\) set to values in \(\{10^{-4}, 2\cdot10^{-4}, 3\cdot10^{-4}\}\).
    \item \emph{Network Pruning:} Post-training, networks were pruned using a threshold \(\varepsilon\) within \(\{10^{-4}, 10^{-3}, 10^{-2}, 10^{-1}\}\) to enhance numerical stability.
\end{enumerate}

Second, for each node in each neural network we proceeded to compute strong and weak activation bounds using Algorithm \ref{alg:bounder} with Gurobi 9.1.2 \cite{gurobi} for MILP and LP runs. A time limit of one hour was used to compute each bound. We also compute the naive activation bounds. 

Third, for each network architecture and each type of activation bounds (strong, weak, naive), we defined instances of the verification problem from Section \ref{sec:verification}. We didn't use multiple seeds, but instead we ran ten instances for each configuration to reduce uncertainty.
For this, we select one network of each architecture with L1 regularization, a pruning threshold $\varepsilon=0.0002$, one random image of each class in MNIST as a reference $x_0$ and a radius $\epsilon=0.18$ to perturb the image. The resulting problem was then solved with Gurobi 9.1.2 \cite{gurobi}. All runs were performed on Dell PowerEdge C6420 Intel Xeon Gold 6152 @ 2.10 GHz with 22 cores; we allocated 4 CPU cores and 6 GB of RAM per network.

Since we would like this study to serve multiple settings, we decided to use a uniform time limit of one hour. For the same reason, we didn't discount the OBBT running time from the verification problem, since its effect depends on the application in mind and how many times the MILP models are solved. Also, the dataset may be a small one, but we hope that advances and insights from smaller problems can help us tackle much larger ones later, as has happened before in other applications of MIPs.

\section{Results}
\subsection{Strong Versus Weak Activation Bounds}

\begin{figure}[t]
    \centering
    \includegraphics[width=0.99\textwidth]{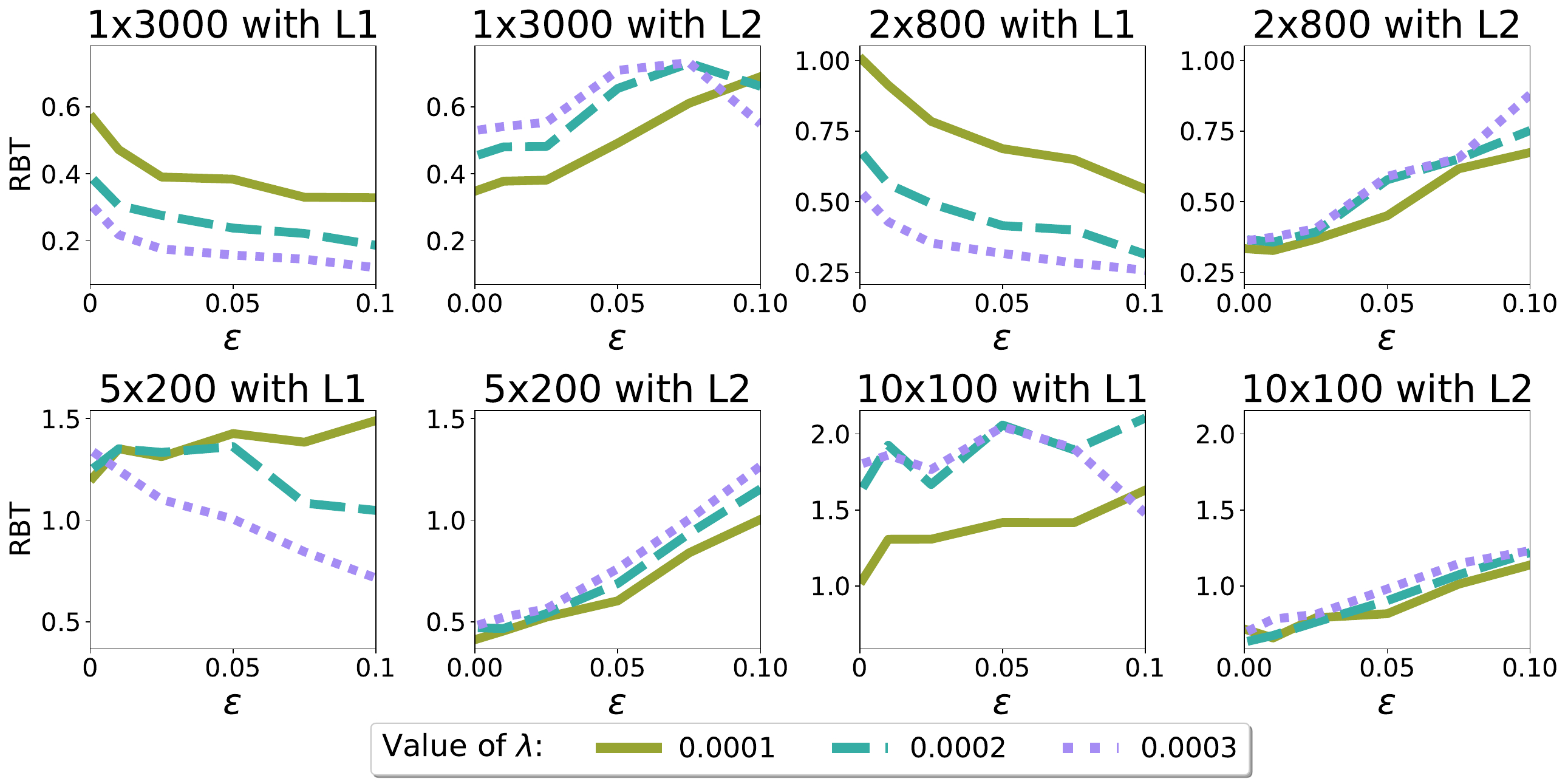}
    \caption{RBT for the lower bound on the last layer versus rounding threshold $\varepsilon$ for each architecture, regularization type (L1/L2) and level $\lambda$.}
    \label{fig:net_quality_reg}
\end{figure}

Our first analysis seeks to determine how different the strong and weak bounds are to each other. This has also been reported in \cite{liu2021algorithms}, but for the sake of completeness we also include our analysis. For this comparison, we use the relative bound tightness (RBT):= $\frac{|\widehat{B} - B|}{|B|+10^{-10}}$, where $\widehat{B}$ and $B$ are the weak and strong activation bounds obtained by the LP and MIP, respectively.
In Figure \ref{fig:net_quality_reg}, we display a summary of our results. Overall, we can see that the strong and weak bounds are remarkably similar to each other. Though not shown, results are similar for the upper activation bounds.
In the worst case, we obtained a weak bound whose value was twice that of the corresponding strong bound for net 10x100.
We also observe an expected decay in the strength of the weak bound with respect to the depth of the network. We confirmed this with an analysis of the RBT layer by layer, which we omit here due to space requirements. However, we remark that we observe an \emph{exponential} growth on RBT by layer.
\begin{figure}[!hbtp]
    \centering
    \includegraphics[width=0.99\textwidth]{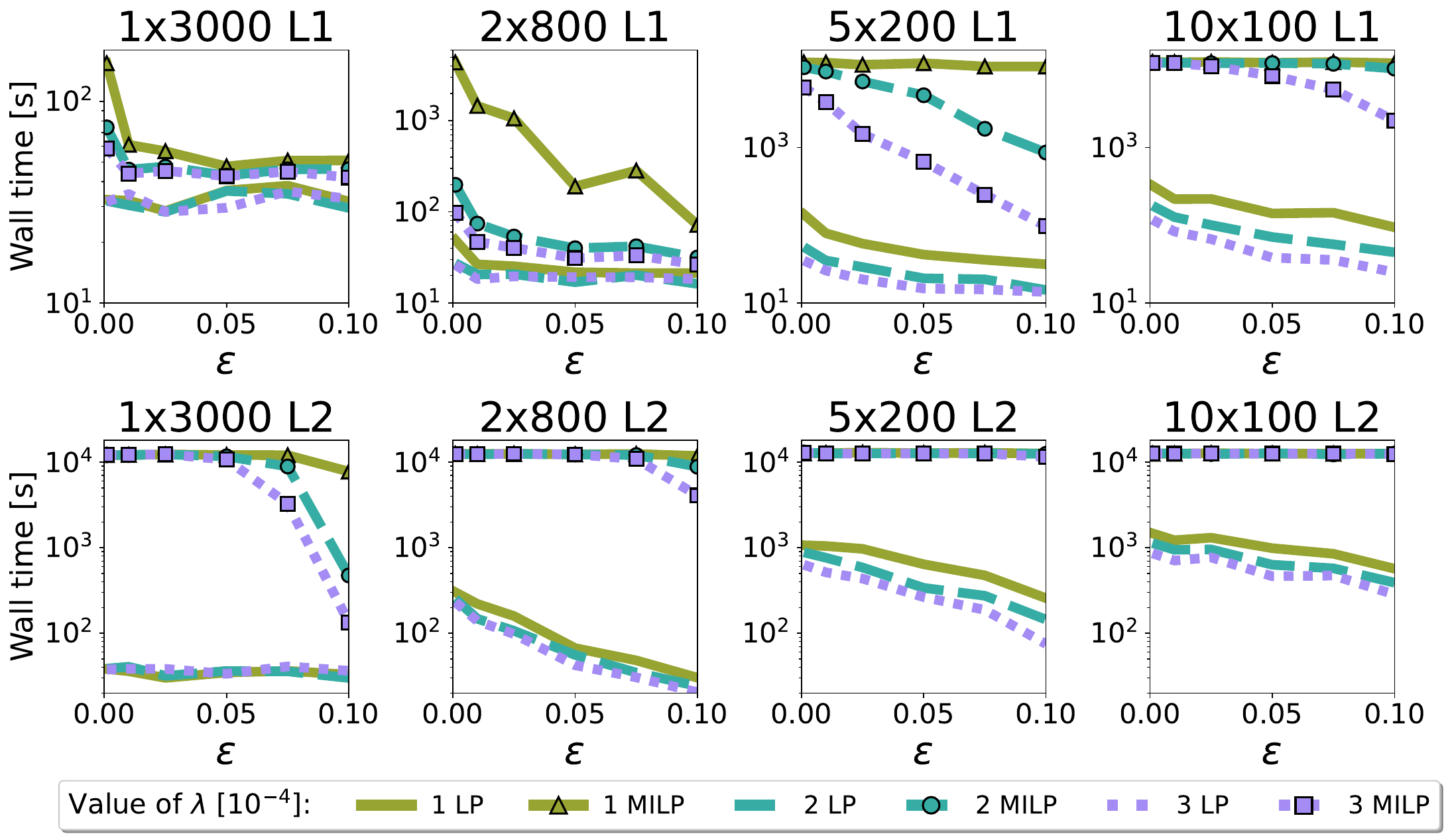}
    \caption{Running times of Algorithm \ref{alg:bounder} for strong and weak bounds in the L1 case (top) and the L2 case (bottom).}
    \label{fig:solv_time_reg}
\end{figure}

\subsection{Solving Times}
A second comparison we report is the running time (wall clock) on the previous experiments. In Figure \ref{fig:solv_time_reg}, we show the total runtime for the L1 and L2 regularization.
In both regularization types, and for every regularization value $\lambda$, we can see the expected improvement in computation time for the weak bounds.
This improvement is much more drastic on deeper networks for the L1 case and more drastic in shallower networks for the L2 case as shown in Figure \ref{fig:solv_time_reg}.
Also note that the pruning threshold has more impact on the strong bounds than on the weak ones, as shown in Figure \ref{fig:solv_time_reg} where the curve drops faster for MILPs.
The same phenomenon happens with $\lambda$ in this case.

Overall, we observe that the weak bounds provide a good alternative for computing activation bounds. They can be computed efficiently, their running time scales well, and they have good quality with respect to the strong bounds. This was also noted in \cite{liu2021algorithms}. Besides confirming this, we also show that the weak bounds performance is robust across different values of pruning thresholds and regularization.
In the following subsection, we test how a worse bound affects a verification model that depends on them. 

\subsection{Effectiveness On Verification Problems}
To test the performance of different bound qualities, we consider three bound types to use in the verification problem from Section \ref{sec:verification}: the strong and weak bounds we discussed before, and \textit{naive} bounds given by Proposition \ref{prop:naive}.

\begin{figure}[t]
    \centering
    \begin{subfigure}[c]{0.4\textwidth}  
        \centering 
        \includegraphics[width=\textwidth]{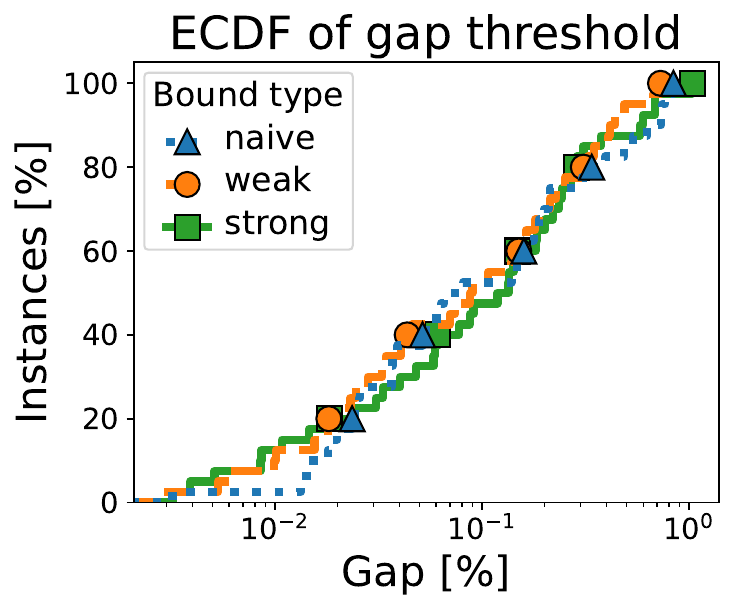}
        \caption{}
        \label{fig:gap_perc}
    \end{subfigure}
    \begin{subfigure}[c]{0.4\textwidth}  
        \centering 
        \includegraphics[width=\textwidth]{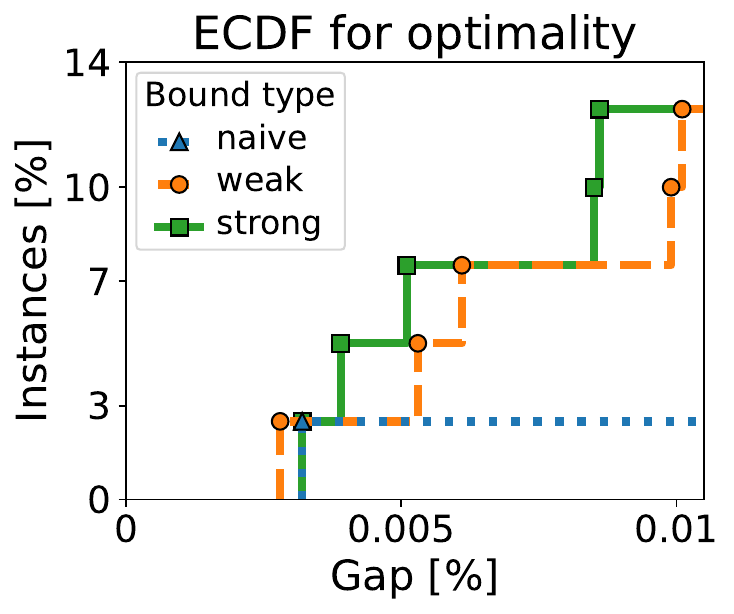}
        \caption{}
        \label{fig:optimal_inst}
    \end{subfigure}
    \caption{Comparison for the different bound types used.}
    \label{fig:bound_type_comp}
\end{figure}
In Figure \ref{fig:optimal_inst}, we see that strong bounds got the most optimal solutions, followed by weak bounds and naive bounds. This was expected. What is more interesting is that if we take a closer look at the instances that were near optimality (gap lower than $1\%$), as shown in Figure \ref{fig:gap_perc}, we can see that weak bounds have a better performance overall, even surpassing strong bounds in some ranges.

\begin{figure}[t]
    \centering 
    \includegraphics[width=\textwidth]{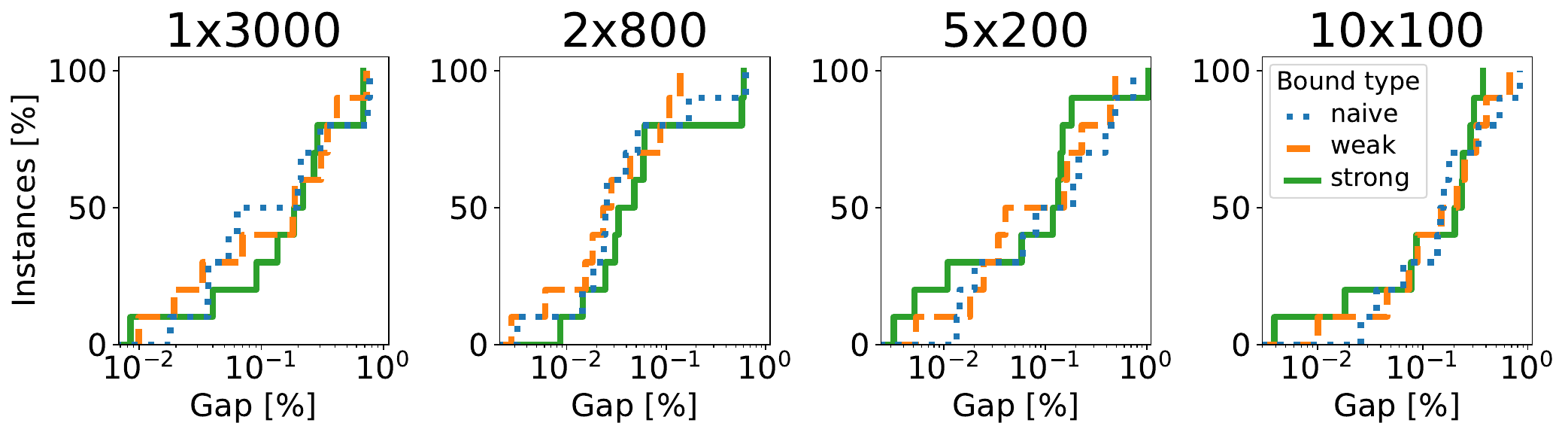}
    \caption{Instances that reach a gap $1\%$ by network architecture.}
    \label{fig:sep_acc}
\end{figure}

In Figure \ref{fig:sep_acc}, we disaggregate the result based on the network architecture. Note that all the bound types reach a gap lower than 1\% for shallower networks. We also note that naive bounds work quite well for these networks. For networks such as 2x800 in Figure \ref{fig:sep_acc}, we can see that weak bounds surpass strong ones, despite taking much lower time to compute. Finally, we see that strong bounds give the best results for the deepest networks overall.
\section{Conclusions}

In terms of the computation of activation bounds, we see that having some properties such as a regularized or pruned network will facilitate their computation. In addition, the weak bounds proved to be comparable to strong bounds, but with considerably less computing time and more stability across different architectures.
Computing the bounds with LP relaxations provides a good compromise of quality in favor of computability.

On the performance of different bound types on a verification problem, the results show that using the weak bounds provides better results overall, and could be used as a general ``rule of thumb''. However, here the dominance is not as clear across different contexts. For example, on shorter networks naive bounds work as well as the other bound types, and on deeper networks the performance of the weak bounds tends to get worse relative to strong bounds.
These results could be used to design a better approach; for example, using a hybrid method that uses naive bounds on the first layers, then weak bounds on the middle layers and strong bounds for the deeper layers.

\begin{credits}
\subsubsection{\ackname} 
G.M. was supported by ANID Fondecyt grant 1231522 and by ANID PIA/PUENTE AFB230002.  M.G. was supported by ANID Fondecyt grant 1231092 and Centro de Modelamiento Matemático (CMM), FB210005, BASAL funds for
centers of excellence from ANID-Chile.
Powered@NLHPC: This research was partially supported by the supercomputing infrastructure of the NLHPC (CCSS210001). 
T.S. was supported by the National Science Foundation (NSF) grant IIS 2104583. 
\end{credits}
%
%
%

%
% ---- Bibliography ----
%
\newpage
\bibliographystyle{splncs04}
\bibliography{bibliography}

\begin{thebibliography}{10}
\providecommand{\url}[1]{\texttt{#1}}
\providecommand{\urlprefix}{URL }
\providecommand{\doi}[1]{https://doi.org/#1}

\bibitem{anderson2018strong}
Anderson, R., Huchette, J., Ma, W., Tjandraatmadja, C., Vielma, J.P.: Strong mixed-integer programming formulations for trained neural networks. Mathematical Programming  \textbf{183}(1-2),  3--39 (2020)

\bibitem{anderson2020strong}
Anderson, R., Huchette, J., Ma, W., Tjandraatmadja, C., Vielma, J.P.: Strong mixed-integer programming formulations for trained neural networks. Mathematical Programming  (2020)

\bibitem{arora2018understanding}
Arora, R., Basu, A., Mianjy, P., Mukherjee, A.: Understanding deep neural networks with rectified linear units. In: International Conference on Learning Representations (ICLR) (2018)

\bibitem{bergman2022janos}
Bergman, D., Huang, T., Brooks, P., Lodi, A., Raghunathan, A.U.: {JANOS}: An integrated predictive and prescriptive modeling framework. INFORMS Journal on Computing  (2022)

\bibitem{burtea2023safe}
Burtea, R.A., Tsay, C.: Safe deployment of reinforcement learning using deterministic optimization over neural networks. In: Computer Aided Chemical Engineering, vol.~52, pp. 1643--1648. Elsevier (2023)

\bibitem{Cacciola2023Jul}
Cacciola, M., Frangioni, A., Lodi, A.: {Structured Pruning of Neural Networks for Constraints Learning}. arXiv  (Jul 2023). \doi{10.48550/arXiv.2307.07457}

\bibitem{cai2023pruning}
Cai, J., Nguyen, K.N., Shrestha, N., Good, A., Tu, R., Yu, X., Zhe, S., Serra, T.: Getting away with more network pruning: From sparsity to geometry and linear regions. In: International Conference on the Integration of Constraint Programming, Artificial Intelligence, and Operations Research (CPAIOR) (2023)

\bibitem{ceccon2022omlt}
Ceccon, F., Jalving, J., Haddad, J., Thebelt, A., Tsay, C., Laird, C.D., Misener, R.: Omlt: Optimization \& machine learning toolkit. Journal of Machine Learning Research  \textbf{23}(349), ~1--8 (2022)

\bibitem{chen2020voltage}
Chen, Y., Shi, Y., Zhang, B.: Data-driven optimal voltage regulation using input convex neural networks. Electric Power Systems Research  (2020)

\bibitem{cheng2017maximum}
Cheng, C.H., N{\"u}hrenberg, G., Ruess, H.: Maximum resilience of artificial neural networks. In: International Symposium on Automated Technology for Verification and Analysis. pp. 251--268. Springer (2017)

\bibitem{cybenko1989approximation}
Cybenko, G.: Approximation by superpositions of a sigmoidal function. Mathematics of Control, Signals and Systems  (1989)

\bibitem{delarue2020rlvrp}
Delarue, A., Anderson, R., Tjandraatmadja, C.: Reinforcement learning with combinatorial actions: An application to vehicle routing. In: NeurIPS (2020)

\bibitem{mnist}
Deng, L.: The mnist database of handwritten digit images for machine learning research. IEEE Signal Processing Magazine  \textbf{29}(6),  141--142 (2012)

\bibitem{elaraby2023oamip}
ElAraby, M., Wolf, G., Carvalho, M.: {OAMIP}: Optimizing {ANN} architectures using mixed-integer programming. In: International Conference on the Integration of Constraint Programming, Artificial Intelligence, and Operations Research (CPAIOR) (2023)

\bibitem{fajemisin2023ocl}
Fajemisin, A., Maragno, D., den Hertog, D.: Optimization with constraint learning: A framework and survey. European Journal of Operational Research  (2023)

\bibitem{fischetti2018deep}
Fischetti, M., Jo, J.: {Deep neural networks and mixed integer linear optimization}. Constraints  \textbf{23}(3),  296--309 (Jul 2018). \doi{10.1007/s10601-018-9285-6}

\bibitem{funahashi1989approximate}
Funahashi, K.I.: On the approximate realization of continuous mappings by neural networks. Neural Networks  (1989)

\bibitem{glorot2011rectifier}
Glorot, X., Bordes, A., Bengio, Y.: Deep sparse rectifier neural networks. In: International Conference on Artificial Intelligence and Statistics (AISTATS) (2011)

\bibitem{gurobi2023ml}
Gurobi: {Gurobi Machine Learning}. \url{https://github.com/Gurobi/gurobi-machinelearning} (2023), accessed: 2023-12-03

\bibitem{gurobi}
{Gurobi Optimization, LLC}: {Gurobi Optimizer Reference Manual} (2023), \url{https://www.gurobi.com}

\bibitem{hahnloser2000origin}
Hahnloser, R., Sarpeshkar, R., Mahowald, M., Douglas, R., Seung, S.: Digital selection and analogue amplification coexist in a cortex-inspired silicon circuit. Nature  (2000)

\bibitem{hanin2017approximating}
Hanin, B., Sellke, M.: Approximating continuous functions by {ReLU} nets of minimal width. arXiv:1710.11278  (2017)

\bibitem{hornik1989approximator}
Hornik, K., Stinchcombe, M., White, H.: Multilayer feedforward networks are universal approximators. Neural Networks  (1989)

\bibitem{huchette2023deep}
Huchette, J., Mu{\ifmmode\tilde{n}\else\~{n}\fi}oz, G., Serra, T., Tsay, C.: {When Deep Learning Meets Polyhedral Theory: A Survey}. arXiv  (Apr 2023). \doi{10.48550/arXiv.2305.00241}

\bibitem{kanamori2021counterfactual}
Kanamori, K., Takagi, T., Kobayashi, K., Ike, Y., Uemura, K., Arimura, H.: Ordered counterfactual explanation by mixed-integer linear optimization. In: AAAI Conference on Artificial Intelligence (AAAI) (2021)

\bibitem{katz2017reluplex}
Katz, G., Barrett, C., Dill, D.L., Julian, K., Kochenderfer, M.J.: Reluplex: An efficient smt solver for verifying deep neural networks. In: Computer Aided Verification: 29th International Conference, CAV 2017, Heidelberg, Germany, July 24-28, 2017, Proceedings, Part I 30. pp. 97--117. Springer (2017)

\bibitem{khalil2018combinatorial}
Khalil, E.B., Gupta, A., Dilkina, B.: Combinatorial attacks on binarized neural networks. In: International Conference on Learning Representations (ICLR) (2019)

\bibitem{lecun2015nature}
LeCun, Y., Bengio, Y., Hinton, G.: Deep learning. Nature  (2015)

\bibitem{li2022sok}
Li, L., Xie, T., Li, B.: Sok: Certified robustness for deep neural networks. In: 2023 IEEE Symposium on Security and Privacy (SP). pp. 94--115. IEEE Computer Society (2022)

\bibitem{liu2021algorithms}
Liu, C., Arnon, T., Lazarus, C., Strong, C., Barrett, C., Kochenderfer, M.J., et~al.: Algorithms for verifying deep neural networks. Foundations and Trends{\textregistered} in Optimization  \textbf{4}(3-4),  244--404 (2021)

\bibitem{lu2017expressive}
Lu, Z., Pu, H., Wang, F., Hu, Z., Wang, L.: The expressive power of neural networks: A view from the width. In: Neural Information Processing Systems (NeurIPS) (2017)

\bibitem{maragno2023mixed}
Maragno, D., Wiberg, H., Bertsimas, D., Birbil, S.I., Hertog, D.d., Fajemisin, A.: Mixed-integer optimization with constraint learning. Operations Research  (2023)

\bibitem{murzakhanov2022powerflow}
Murzakhanov, I., Venzke, A., Misyris, G.S., Chatzivasileiadis, S.: Neural networks for encoding dynamic security-constrained optimal power flow. In: Bulk Power Systems Dynamics and Control Sympositum (2022)

\bibitem{nair2010rectified}
Nair, V., Hinton, G.: Rectified linear units improve restricted boltzmann machines. In: International Conference on Machine Learning (ICML) (2010)

\bibitem{park2021width}
Park, S., Yun, C., Lee, J., Shin, J.: Minimum width for universal approximation. In: International Conference on Learning Representations (ICLR) (2021)

\bibitem{pascanu2013on}
Pascanu, R., Mont\'{u}far, G., Bengio, Y.: On the number of response regions of deep feedforward networks with piecewise linear activations. In: International Conference on Learning Representations (ICLR) (2014)

\bibitem{perakis2022optimizing}
Perakis, G., Tsiourvas, A.: Optimizing objective functions from trained relu neural networks via sampling. arXiv:2205.14189  (2022)

\bibitem{ramachandran2018pop}
Ramachandran, P., Zoph, B., Le, Q.V.: Searching for activation functions. In: ICLR Workshop Track (2018)

\bibitem{say2017planning}
Say, B., Wu, G., Zhou, Y.Q., Sanner, S.: Nonlinear hybrid planning with deep net learned transition models and mixed-integer linear programming. In: International Joint Conference on Artificial Intelligence (IJCAI) (2017)

\bibitem{serra2020empirical}
Serra, T., Ramalingam, S.: Empirical bounds on linear regions of deep rectifier networks. In: AAAI Conference on Artificial Intelligence (AAAI) (2020)

\bibitem{serra2020lossless}
Serra, T., Kumar, A., Ramalingam, S.: Lossless compression of deep neural networks. In: 17th International Conference on Integration of Constraint Programming, Artificial Intelligence, and Operations Research (CPAIOR). pp. 417--430. Springer (2020)

\bibitem{serra2017bounding}
Serra, T., Tjandraatmadja, C., Ramalingam, S.: Bounding and counting linear regions of deep neural networks. In: International Conference on Machine Learning (ICML). pp. 4558--4566. PMLR (2018)

\bibitem{serra2021scaling}
Serra, T., Yu, X., Kumar, A., Ramalingam, S.: Scaling up exact neural network compression by {ReLU} stability. Advances in neural information processing systems  \textbf{34},  27081--27093 (2021)

\bibitem{Tjeng2017Nov}
Tjeng, V., Xiao, K., Tedrake, R.: {Evaluating Robustness of Neural Networks with Mixed Integer Programming}. arXiv  (nov 2017). \doi{10.48550/arXiv.1711.07356}

\bibitem{tong2023walk}
Tong, J., Cai, J., Serra, T.: Optimization over trained neural networks: Taking a relaxing walk  (2023)

\bibitem{tsay2021partition}
Tsay, C., Kronqvist, J., Thebelt, A., Misener, R.: Partition-based formulations for mixed-integer optimization of trained {ReLU} neural networks. In: Neural Information Processing Systems (NeurIPS). vol.~34 (2021)

\bibitem{wu2020scalable}
Wu, G., Say, B., Sanner, S.: Scalable planning with deep neural network learned transition models. Journal of Artificial Intelligence Research  (2020)

\bibitem{yang2021control}
Yang, S., Bequette, B.W.: Optimization-based control using input convex neural networks. Computers \& Chemical Engineering  (2021)

\bibitem{yarotsky2017relu}
Yarotsky, D.: Error bounds for approximations with deep {ReLU} networks. Neural Networks  (2017)

\end{thebibliography}

\end{document}